\documentclass{amsart}
\usepackage[english,french]{babel}
\usepackage{fullpage}

\usepackage{lmodern}
\usepackage{amsfonts}
\usepackage[T1]{fontenc}

\usepackage[latin1]{inputenc}

\usepackage{epsfig}
\usepackage{color,graphicx}
\usepackage{amsmath}
\usepackage{amssymb}
\usepackage{varioref}
\usepackage{maththm}

\usepackage[bookmarks=false]{hyperref}

\def\doi#1{\textsl{\href{http://dx.doi.org/#1}{DOI: #1}}}

\newcommand{\ud}{\,\mathrm{d}}

\newcommand{\R}{\mathbb{R}}

\newcommand{\on}{\operatorname}
\newcommand{\real}{\mathbb{R}}

\newcommand{\ip}[3]{\left< {#1}, {#2} \right>_{#3}}

\newcommand{\St}{\mathbf{St}}
\newcommand{\Gr}{\mathbf{Gr}}

\newcommand{\comment}[1]{ }

\makeatletter
\def\len{\mathop{\operator@font len}\nolimits}
\def\Len{\mathop{\operator@font Len}\nolimits}
\def\dil{\mathop{\operator@font dil}\nolimits}
\def\mod{\mathop{\operator@font mod}\nolimits}
\def\supess{\mathop{\operator@font supess}\nolimits}
\makeatother

\begin{document}
\selectlanguage{english}

\title{
Geodesics in infinite dimensional Stiefel and Grassmann manifolds /
Géodesiques sur des variétiés de Stiefel et de Grassmann de dimension infinie
}
\thanks{This research was funded by \texttt{SNS09MENNB} of the Scuola Normale
Superiore, and by the FWF Project 21030.}

\author[P. Harms]{Philipp Harms
\address{P. Harms, Harvard Education Innovation Laboratory}
\email{pharms@edlabs.harvard.edu}
}
\author[A. Mennucci]{Andrea C. G. Mennucci
\address{A. Mennucci, Scuola Normale Superiore,
Piazza dei Cavalieri 7, 56126, Pisa, Italy}
\email{a.mennucci@sns.it}
}

\date{}

\begin{abstract} 
  Let $V$ be a separable Hilbert space, possibly infinite dimensional.
  Let $\St(p,V)$ be the Stiefel manifold of orthonormal frames of
  $p$ vectors in $V$, and let $\Gr(p,V)$ be the Grassmann manifold of $p$
  dimensional subspaces of $V$. We study the distance and the
  geodesics in these manifolds, by reducing the matter to the finite dimensional
  case. We then prove that any two points in those manifolds can be
  connected by a minimal geodesic, and characterize the cut locus.
  
  \medskip
  
  \selectlanguage{french} \noindent \textsc{R\'esum\'e}.  
  Soit $V$ un espace de Hilbert séparable, éventuellement de
  dimension infinie.  Soient $\St(p,V)$ l'ensemble des systèmes
  orthonormés de $p$ vecteurs de $V$, appelé la variété de
  Stiefel, et $\Gr(p,V)$ l'ensemble des sous-espaces vectoriels de $V$
  de dimension $p$, appelé la variété Grassmannienne.  En
  réduisant le problème en dimension finie, nous montrons que dans
  ces espaces il existe des géodésiques minimales entre chaque
  paire de points et nous caractérisons le cut-locus.
  \selectlanguage{english}
\end{abstract}

\maketitle

\section{Introduction}

Let $V$ be a separable Hilbert space, let $p$ be a positive natural number.
We assume that $\dim(V)\ge (2p)$ from here on.
$\St(p,V)$ is the set of orthonormal frames of $p$ vectors in $V$.
Equivalently, we consider 
$$\St(p,V)=\{ x \in L(\R^p,V): x^\top \circ x = \on{Id}_{\R^p} \}$$ 
to be the set of all linear isometric immersions of $\real^p$ into $V$.
Here $x^\top \in L(V,\R^p)$ is the transpose with respect to the metrics on $V$ and $\R^p$, i.e.
$$\big\langle x^\top ( v), r \big\rangle_{\R^p} 
	= \big\langle v, x(r) \big\rangle_{V} \text{ for all } v \in V, r \in \R^p.$$
$\St(p,V)$ is a smooth embedded submanifold in $V^p$. 
The induced Riemannian metric on $\St(p,V)$ is $\ip{x}{y}{} =  \mbox{tr}(x^\top y)$. 
$\St(p,V)$ is a complete Riemannian manifold with this metric.
$\Gr(p,V)$ is the manifold of $p$-dimensional linear subspaces of $V$
and equals the orbit space $\St(p,V)/O(p)$ with respect to $O(p)$ acting
on $\St(p,V)$ by composition from the right.

Our interest is due to the fact that 
$\St(2,V)$ with $V=L^2([0,1])$ is isometric to
the space of planar closed curves up to translation and scaling,
endowed with a Sobolev metric of order one. 
The $O(2)$-action on $\St(2,V)$ corresponds to rotations of the curves. 
Thus $\Gr(2,V)$ with $V=L^2([0,1])$ is isometric to the space of planar closed curves 
up to translations, scalings and rotations.
See \cite{Younes:Comp}, \cite{MR2383560}, \cite{SMSY09:cdc} and \cite{SMSY10SIAM}.
Any result that is proven about the Stiefel or Grassmannian
immediately carries over to the corresponding space of curves. 

\section{Critical geodesics}

We will call a curve $\gamma$ in a Riemannian manifold a
\textbf{critical geodesic} if it is a solution to
the equation $\nabla_{\partial_t} {\dot \gamma}=0$, where $\nabla$ is the covariant derivative.
Such $\gamma$ is a critical point for the action
$\int_0^1 g(\dot \gamma,\dot\gamma ) \ud t$.

\begin{Proposition}[Critical geodesics in $\St(p,V)$ ]
Let $\St(p,V)$ be endowed with the induced metric from $V^p$. 
Let $\gamma : [0,1] \to \St(p,V)$ be a path. Then the geodesic equation is
$\ddot{\gamma} + \gamma(\dot{\gamma}^\top\dot{\gamma}) = 0$.
Solutions to the geodesic equation exist for all time and are given by
  \begin{equation} \label{eq:Stiefel_soln_rn}
    (\gamma(t)e^{At}, \dot{\gamma}(t)e^{At}) = (\gamma(0), \dot{\gamma}(0)) 
    \exp{ t
      \left(
      \begin{array}{ll}
        A         & -S \\
        \mbox{Id} &  A
      \end{array}
      \right)
    }
  \end{equation}
  where $A = \gamma(0)^\top\dot{\gamma}(0)$, $S = \dot{\gamma}(0)^\top\dot{\gamma}(0)$, and
  $\mbox{Id}$ is the $p\times p$ identity matrix.
\end{Proposition}
For $V=\R^n$ this has been demonstrated by Edelman et
al.~\cite[section~2.2.2]{EdelmanAriasSmith98}%
\footnote{\cite{EdelmanAriasSmith98} credits
a personal communication by R.~A.~Lippert
for the final closed form formula \eqref{eq:Stiefel_soln_rn}.}. 
Going through their proof one sees that the same result holds when $V$ is infinite dimensional.

\begin{Proposition}\label{prop:span}
  Equation \eqref{eq:Stiefel_soln_rn} shows that the subspace of $V$ spanned 
  by the $(2p)$ columns of $\gamma(t),\dot \gamma(t)$ remains 
  in the space spanned by the columns of
  $\gamma(0),\dot{\gamma}(0)$ for all $t$.

  This means that, if  $W$ is the subspace  of $V$
  spanned by the columns of  $\gamma(0),\dot{\gamma}(0)$, then
  we can formulate the geodesic equation as an equation in 
  $\St(2,W)$. Obviously, $\dim(W)\le 2p$.
  
  This also means that, if $\gamma$ is a critical geodesic connecting
  $x$ to $y$, and the space $W$ spanned by the columns of $x,y$ is
  $(2p)$ dimensional, then, for any $t$, the columns of $\gamma(t)$ and
  of $\dot \gamma(t)$ must be contained in $W$.
\end{Proposition}

\section{Minimal geodesics}

We denote by $d(x,y)$
the infimum of the length of all paths connecting two points $x,y$
in a Riemannian manifold. 
It does not matter whether the infimum is taken over smooth
or absolutely continuous paths~%
\footnote{Lemma 6.1 in Chap.~VIII in \cite{Lang:FDG}
  can be used to convert any absolutely continuous path
  to a shorter and piecewise smooth path.}.
We call a path $\gamma$ a  \textbf{minimal geodesic} if its length
is equal to the distance $d(\gamma(0),\gamma(1))$. 
Up to a time reparametrization,
a minimal geodesic is smooth and is a critical geodesic.
We will always silently assume that minimal geodesics are parametrized 
such that they are critical.

Let $(M,g)$ be a Riemannian manifold, and $d$ be the induced distance.
When $M$ is finite dimensional, by the celebrated
Hopf--Rinow theorem, metric completeness of $(M,d)$
is equivalent to geodesic completeness of $(M,g)$, and both
imply that any two points $x,y\in M$ can be connected by a minimal geodesic.
In infinite dimensional manifolds this is not true in general.
Indeed, in \cite{atkin75} there is an example
of an infinite dimensional
metrically complete Hilbert smooth manifold $M$ and $x,y\in M$ such that there
is no critical (and thus no minimal) geodesic connecting $x$ to $y$. 
A simpler example, due to Grossman \cite{Grossman65}
(see also sec. VIII.\S6 in \cite{Lang:FDG}),
is an infinite dimensional ellipsoid where the south and north
pole can be connected by countably many critical geodesics of decreasing
length, so that the distance between the poles is not attained by any minimal geodesic.

We will show that, even when $V$ is infinite dimensional,
any two points in $\St(p,V)$ and $\Gr(p,V)$ can be connected
by a minimal geodesic.

\subsection{Minimal geodesics in the Stiefel manifold}

\begin{Theorem}\label{thm:geo_St}
	Let $V$ be a Hilbert space.
  Consider a $(2p)$ dimensional Hilbert space $W$ and an isometric linear embedding
  $i:W\to V$. Then $i$ induces an isometric embedding 
  $$i_*:\St(p,W) \to \St(p,V), \qquad x \mapsto i \circ x $$ 
  (here we consider $x\in \St(p,W)$ to be 
  a linear isometric immersion of $\real^p$ into $W$).
  \begin{enumerate}
  \item  $i_*\big(\St(p,W))$ is totally geodesic in $\St(p,V)$.
  
  \item 
    Let $d_W$ be the distance in $\St(p,W)$
    and similarly $d_V$ in $\St(p,V)$, then
    \begin{equation}
      \label{eq:d_2p_d_n}
      d_W(x,y)=d_V\big(i_*(x),i_*(y)\big)~~.
    \end{equation}
      
  \item 
    Let $x,y\in \St(p,W)$, and a minimal geodesic $\gamma$
    connecting $x$ to $y$ in $\St(p,W)$: then $i_* \circ \gamma$
    is a minimal geodesic connecting $i_*( x)$ to $i_*(y)$ in
    $\St(p,V)$.
	
  \item
    The diameter of $\St(p,V)$ is equal to the diameter of $\St(p,\R^{2p})$.

  \item Any two points $x,y \in \St(p,V)$ can be connected by a minimal geodesic
    $\gamma$. 
    Any minimal geodesic lies in $\St(p,U)$, where $U$ is a $(2p)$ dimensional 
    subspace of $V$ (dependent on $\gamma$).
    
  \item Let $x,y\in \St(p,V)$. Then $y$ is in the cut locus of $x$ 
  	if and only if there is a $(2p)$ dimensional subspace $W$ of $V$ and $\tilde x, \tilde y \in \St(p,W)$ 
	such that $x=i_*(\tilde x)$, $y=i_*(\tilde y)$ and
	$i_*(y)$ is in the cut locus of $i_*(x)$.	
  \end{enumerate}
\end{Theorem}  
Note that point \emph{(5)} in 
the above theorem implies that minimal geodesics can
be numerically computed using a finite dimensional algorithm;
see Sec.~3.3.4 in \cite{SMSY10SIAM}.
\smallskip

We will need two lemmas.

\begin{Lemma}\label{lem:st_linearly_independent}
  Given $x\in \St(p,V)$, the set of $y\in \St(p,V)$ such that
  the columns of $x,y$ are linearly independent is dense in
  $\St(p,V)$.
\end{Lemma}
\begin{proof}
  Let $U$ be the linear space spanned by  the columns of $x,y$;
  if this space is not $(2p)$ dimensional, then let $r_1,\ldots r_k$
  be orthonormal vectors that lie in $U^\perp$, with $k=2p-\dim(U)$.
  Up to  reindexing the columns of $y$, we can suppose that
  the columns $x_1,\ldots,x_p,y_1,\ldots,y_{p-k}$ are linearly independent.
  For  $\varepsilon>0$ small, we then define
  \[\tilde y_i=
  \begin{cases}
    y_i & i=1,\ldots p-k \\
    \cos(\varepsilon) y_i+ \sin (\varepsilon) r_i & i=(p-k+1),\ldots p
  \end{cases}~~.
  \]

  It is easy to verify that $\tilde y\in  \St(p,V)$ and that
  the columns of $x,\tilde y$ are linearly independent.
\end{proof}

\begin{Lemma}\label{lem:st_finite}
  The theorem holds when $V$ is a Hilbert space of finite dimension
  $n$ with $n>2p$.
\end{Lemma}

\begin{proof}
  We prove point \emph{(1)}.
    Let us consider the subgroup $G=O\big(i(W)^\perp\big)$ of $O(V)$
    that keeps $i(W)$ fixed. Then $G$ acts isometrically on
    $\St(p,V)$ as well, and its fixed point set is
    $i_*\big(\St(p,W)\big)$. This proves that $i_*\big( \St(p,W)\big)$
    is totally geodesic in $\St(p,V)$.
  
	To prove point \emph{(2)} we first note that 
	since $\St(p,W)$ is isometrically embedded in $\St(p,V)$, we have
	$d_W(x,y) \geq d_V\big(i_*(x),i_*(y)\big)$.
    We will show the inverse inequality only for the case when the columns of 
    $x$ and $y$ are linearly independent.
    The general case then follows because the set of $y$ such that the columns of $x$ and $y$ are linearly
    independent is dense in $W$ by lemma~\ref{lem:st_linearly_independent} and
    since distances are Lipschitz continuous. 

    Since $V$ is finite dimensional, $\St(p,V)$ is compact, so 
    by the Hopf--Rinow Theorem
    $i_*(x)$ and $i_*(y)$ can be connected by a minimizing geodesic in
    $\St(p,V)$. The columns of $i_*(x)$ and $i_*(y)$ together span 
    the $(2p)$ dimensional space $i(W)$, so we can apply
    proposition \ref{prop:span}. This allows us to write $\gamma=i_* \circ \tilde \gamma$ for
    a path $\tilde \gamma$ in $\St(p,W)$. Then
    $$d_W(x,y)\leq \len(\tilde \gamma) = \len(i_* \circ \tilde \gamma) = \len(\gamma)= 
       d_V\big(i_*(x),i_*(y)\big) ~~.$$

    Point \emph{(3)} follows from point \emph{(2)} and the equality
    \[\len(i_*\circ\gamma)=\len(\gamma)=d_W(x,y)=d_V\big(i_*(x),i_*y\big)~~.\]

    Point \emph{(4)} follows from point \emph{(3)}. 
    Point \emph{(5)} follows from the Hopf--Rinow theorem and 
    the discussion in Prop.~\ref{prop:span}.

  	We now prove point \emph{(6)}.
	By definition, $y$ is in the cut locus of $x$ if and only if there is 
	a geodesic $\gamma$ in $\St(p,V)$ with $\gamma(0)=x, \gamma(1)=y$ such that
	$$\sup \big\{ t: \len(\gamma|_{[0,t]}) = d_V\big(\gamma(0),\gamma(t) \big) \big\} = 1.$$
	(Recall that we write $d_V$ for the distance in $\St(p,V)$.)
	Any such geodesic lies in $\St(p,W)$ for some $(2p)$ dimensional space $W$. 
	Letting $i:W \to V$ denote the isometric embedding, we can write 
	$\gamma = i_* \circ \tilde \gamma$ for a path $\tilde \gamma$ in $\St(p,W)$.
	Then one has by point \emph{(2)} that
	\begin{equation*}
		\sup \big\{ t: \len(\tilde \gamma|_{[0,t]}) = d_W\big(\tilde\gamma(0),\tilde\gamma(t) \big) \big\} = 1~~.		\qedhere
	\end{equation*}
\end{proof}

We now prove Theorem~\ref{thm:geo_St}.

\begin{proof}
  The proof of points \emph{(1), (3), (4), (6)} works as in the finite dimensional case. 
  We will now prove point \emph{(2)}. We have
  $d_W(x,y) \geq d_V\big(i_*(x),i_*(y)\big),$
  since $\St(p,W)$ is isometrically embedded in $\St(p,V)$.
  It remains to show the inverse inequality. 

  Consider a smooth path $\xi$ connecting $i_*(x)$ to $i_*(y)$ in $\St(p,V)$.
  We can find finitely many points $0=t_0<t_1<\ldots <t_k=1$
  such that $\xi|_{{[t_i,t_{i+1}]}}$ is contained 
  inside (the manifold part of) a  normal  neighborhood.
  So $\xi(t_i),\xi(t_{i+1})$ can be connected by a minimal geodesic. 
  By joining all these minimal geodesics we obtain a piecewise smooth path 
  $\eta$, with $\len(\eta)\le\len(\xi)$.
  Then by repeated application of proposition~\ref{prop:span} there is a finite 
  dimensional subspace $\tilde W$ of $V$ that contains the columns of $\eta(t)$ for $t\in[0,1]$. 
  When necessary we enlarge $\tilde W$ such that it also contains $i(W)$.
  Now the finite dimensional version of this Lemma allows 
  us to compare $\St(p,W)$ to $\St(p,\tilde W)$, and we get: 
  $$d_W(x,y) = d_{\tilde W}\big(i_*(x),i_*(y)\big) \leq \len(\eta) \leq \len(\xi)~~.$$
  Since this holds for arbitrary 
  paths $\xi$ connecting $i_*(x)$ to $i_*(y)$ in $\St(p,V)$, 
  we get 
  $d_W(x,y) \leq d_V\big(i_*(x),i_*(y)\big)$.
		
  Point \emph{(5)}
  now follows by choosing any linear subspace $W$ containing the columns of $x,y$.
\end{proof}

\subsection{Minimal geodesics in the Grassmann manifold}

\begin{Theorem}\label{thm:geo_Gra}
	Thm.~\ref{thm:geo_St} remains valid when Stiefels are replaced by Grassmannians.
	Most importantly, for any two points $x,y \in\Gr(p,V)$, 
	there is a minimal geodesic $\gamma$ connecting $x$ to $y$.
	The same holds for the  Grassmannian $\Gr_+(p,V)$ of oriented $p$ spaces.
\end{Theorem}

We need a Lemma.  
\begin{Lemma}[Existence of horizontal paths]\label{lem:horizontal_paths}
	For any path $x:[0,1]\rightarrow \St(p,V)$ there is a path $g:[0,1] \rightarrow O(p)$
	such that the path $x(t) \circ g(t)$ is horizontal, i.e. normal 
	to the $O(p)$-orbits in $\St(p,V)$.
\end{Lemma}

\begin{proof}
	The tangent space at $x$ to the $O(p)$-orbit through $x \in \St(p,V)$ is 
	$$T_x\big(x . O(p)\big) 
		= \big\{ x z : z \in \mathfrak o(p) \big\}
		= \big\{ x z : z \in L(\R^p,\R^p), z^\top + z = 0 \big\}~~.$$
	Thus $y \in T_x\St(p,V)$ is horizontal if and only if 
	$\mbox{tr}(y^\top x z)=0$ for all antisymmetric $z$. 
	This is equivalent to $y^\top x$=0 
	because $y^\top x$ is antisymmetric, too. 
	Thus the path $x g$ is horizontal iff
	$$\big(\partial_t ( x g )\big)^\top (x g)
		= g^\top\dot x^\top x g+\dot g^\top x^\top x g =0.$$
	This can be achieved by letting $g$ be the solution to the ODE
	$\dot g = - x^\top \dot x g$.
\end{proof}
Note that the length of $x(t) \circ g(t)$ is smaller than or equal to the length of 
$x(t)$, with equality if and only if $x(t)$ is already a horizontal path.

We are now able to prove Thm.~\ref{thm:geo_Gra}. 
\begin{proof}
$\St(p,V)$ is a principal fiber bundle with structure group $O(p)$
over $\Gr(p,V)=\St(p,V)/O(p)$.
We prove the existence of minimizing geodesics connecting any two 
points in $\Gr(p,V)$.
Take any point $\tilde x \in \St(p,V)$
in the fiber over $x$. The fiber over $y$ is compact since $O(p)$ is compact. 
Therefore $d(x,\cdot)$ attains a minimum at some point $\tilde y$ in the fiber 
over $y$. By theorem~\ref{thm:geo_St} there is a minimal geodesic connecting 
$\tilde x$ to $\tilde y$. This geodesic is horizontal since otherwise it could
be made shorter by making it horizontal. (We use lemma~\ref{lem:horizontal_paths} here.) 
By the theory of Riemannian submersions it projects to a minimal geodesic in $\Gr(p,V)$.

The remaining statements simply follow from Thm.~\ref{thm:geo_St} by
going to the quotient with respect to the $O(p)$-action.
For the case of  $\Gr_+(p,V)$, we use the group $SO(p)$ instead of $O(p)$.
\end{proof}

\bibliographystyle{myplain}
\bibliography{mennbib} 

\begin{thebibliography}{1}

\bibitem{atkin75}
C.~J. Atkin.
\newblock The {H}opf-{R}inow theorem is false in infinite dimensions.
\newblock {\em Bull. London Math. Soc.}, 7(3):261--266, 1975.
\newblock \doi{10.1112/blms/7.3.261}.

\bibitem{EdelmanAriasSmith98}
A.~Edelman, T.~Arias, and S.~Smith.
\newblock The geometry of algorithms with orthogonality constraints.
\newblock {\em SIAM J Matrix Analy Appl}, 20:303--353, 1998.
\newblock \doi{10.1137/S0895479895290954}.
\newblock arXiv:physics/9806030v1 (1998).

\bibitem{Grossman65}
Nathaniel Grossman.
\newblock Hilbert manifolds without epiconjugate points.
\newblock {\em Proc. Amer. Math. Soc.}, 16:1365--1371, 1965.
\newblock ISSN 0002-9939.
\newblock \doi{10.2307/2035933}.

\bibitem{Lang:FDG}
Serge Lang.
\newblock {\em Fundamentals of differential geometry}.
\newblock Springer--Verlag, 1999.
\newblock ISBN 0-387-98593-X.

\bibitem{SMSY09:cdc}
Ganesh Sundaramoorthi, Andrea Mennucci, Stefano Soatto, and Anthony Yezzi.
\newblock Tracking deforming objects by filtering and prediction in the space
  of curves.
\newblock In {\em Conference on Decision and Control}, pages 2395 -- 2401,
  2009.
\newblock ISBN 978-1-4244-3871-6.
\newblock \doi{10.1109/CDC.2009.5400786}.

\bibitem{SMSY10SIAM}
Ganesh Sundaramoorthi, Andrea Mennucci, Stefano Soatto, and Anthony Yezzi.
\newblock A new geometric metric in the space of curves, and applications to
  tracking deforming objects by prediction and filtering.
\newblock {\em SIAM Journal on Imaging Sciences}, 4:109--145, 2011.
\newblock \doi{10.1137/090781139}.

\bibitem{Younes:Comp}
Laurent Younes.
\newblock Computable elastic distances between shapes.
\newblock {\em SIAM Journal of Applied Mathematics}, 58(2):565--586, 1998.
\newblock \doi{10.1137/S0036139995287685}.

\bibitem{MR2383560}
Laurent Younes, Peter~W. Michor, Jayant Shah, and David Mumford.
\newblock A metric on shape space with explicit geodesics.
\newblock {\em Atti Accad. Naz. Lincei Cl. Sci. Fis. Mat. Natur. Rend. Lincei
  (9) Mat. Appl.}, 19(1):25--57, 2008.
\newblock ISSN 1120-6330.
\newblock \doi{10.4171/RLM/506}.

\end{thebibliography}

\end{document}